\documentclass[fleqn]{mat01}
\usepackage{times,mathtimy,amssymb,latexsym,amscd}
\begin{document}

\setcounter{page}{51}
\firstpage{51}

\newtheorem{thm}{\bf Theorem}[ssection]

\title{The Jacobian of a nonorientable Klein surface, II}

\markboth{Pablo Ar\'es-Gastesi and Indranil Biswas}{Jacobian of a Klein surface, II}

\author{PABLO AR\'ES-GASTESI and INDRANIL BISWAS}

\address{School of Mathematics, Tata Institute of 
Fundamental Research, Mumbai~400~005, India\\
\noindent E-mail: pablo@math.tifr.res.in; indranil@math.tifr.res.in}

\volume{116}

\mon{February}

\parts{1}

\pubyear{2006}

\Date{MS received 23 August 2005}

\begin{abstract}
The aim here is to continue the investigation
in \cite{AB} of Jacobians of a Klein surface
and also to correct an error in \cite{AB}.
\end{abstract}

\keyword{Klein surface; divisor; Jacobian variety.}

\maketitle

\section{Introduction}

This note is a continuation of the study of the Jacobian of
nonorientable, compact Klein surfaces that we started in \cite{AB}. To
explain our results we need to recall the set-up in \cite{AB}. Let $Y$
be a nonorientable, connected compact Klein surface (the changes of
coordinates in $Y$ are either holomorphic or anti-holomorphic), and let
$X$ be the canonical double cover of $Y$ given by local orientations. It
is well-known that $X$ is a Riemann surface. The nontrivial deck
transformation, which we will denote by $\sigma$, for covering $X
\rightarrow Y$ is an anti-holomorphic involution on $X$ with $Y$
as the quotient. 

The Jacobian of $X$ was considered in \cite{AB} from three different
points of view: (1)~as the dual of the space of holomorphic $1$-forms,
(2)~as the divisor class group, and (3)~as the group of isomorphism
classes of degree zero line bundles. In each case, the involution
$\sigma$ induces an involution in the corresponding description of the
Jacobian. In Theorem 4.2 of \cite{AB} it was stated that the Jacobian of
$Y$ can be identified with the set of line bundles $L$ such that
$\sigma^*\overline{L}$ is holomorphically isomorphic to $L$. However,
the proof of that result is not correct. It was asserted in the proof
that the isomorphism of $L$ given by $\sigma^*\overline\alpha \circ
\alpha$ can be so chosen that it is an involution of the total space of
$L$. In this note we show that the isomorphism can be chosen to be of
order four, and there may not be any isomorphism of order two. 

Let $\mathcal G$ denote the group of line bundles $L$ over $X$ such that
$\sigma^*\overline{L}$ is holomorphically isomorphic to $L$. Let
${\mathcal G}_0\, \subset\, {\mathcal G}$ be the subgroup defined by all
$L$ such that the isomorphism of $L$ with $\sigma^*\overline{L}$ can be
chosen to be an involution of the total space of $L$. 

The correct version of Theorem 4.2 of \cite{AB} is as follows:

\setcounter{section}{2}
\setcounter{thm}{4}
\begin{thm}[\!]
The image of the homomorphism $\rho$ in p.~$147${\rm ,} eq.~{\rm
(}$3${\rm )} of {\rm \cite{AB}} coincides with $J^0(X)\cap
{\mathcal G}_0$. If the genus of $X$ is even{\rm ,} then the image of
the homomorphism $\rho$ {\rm (}p.~$147${\rm ,} eq.~{\rm (}$3${\rm )} of
{\rm \cite{AB}}{\rm )} coincides with $J^0(X)\cap {\mathcal G}$. If
the genus of $X$ is odd{\rm ,} then the image of $\rho$ is a subgroup of
$J^0(X)\cap {\mathcal G}$ of index two.
\end{thm}

This theorem is proved in \S~1. In \S\S~2 and 3 we
considered two other approaches to the Jacobian mentioned above (dual
space to the space of holomorphic forms and divisors) and proved this
result in those settings. Thus we get three different (but clearly
equivalent) proofs of the theorem, which shows that one can consider any
of these three ways of describing the Jacobian of $X$ to study the
Jacobian of $Y$.

\setcounter{section}{1}
\section{Real points of the Jacobian and real algebraic line
bundles}

Let $X$ be a compact connected Riemann surface
equipped with a fixed-point free
anti-holomorphic involution $\sigma$. So
$\sigma\hbox{:}\ X \rightarrow X$ is a
diffeomorphism with the property that if
$J(X)\, \in\, C^\infty(X, \text{End}(T^{\mathbb R}X))$
is the almost complex structure on $X$, then
$\sigma^*J(X) = -J(X)$. Note that this implies that
$\sigma$ is orientation reversing.

Let $L$ be a holomorphic line bundle over
$X$ such that $\sigma^*\overline{L}$ is holomorphically
isomorphic to $L$. We note that $\sigma^*\overline{L}$
is a holomorphic line bundle (see \S~4 of \cite{AB}
for the details).

Take an isomorphism
\begin{equation}\label{s.c.1}
\alpha\hbox{:}\ L \rightarrow \sigma^*\overline{L}\,.
\end{equation}
Therefore,
\begin{equation}\label{s.c.2}
(\sigma^*\overline{\alpha})\circ \alpha\hbox{:}\ L \rightarrow L
\end{equation}
is an automorphism of the holomorphic line bundle $L$. Indeed, we have
\begin{equation*}
\sigma^*\overline{\alpha}\hbox{:}\ \sigma^*\overline{L}\rightarrow
\sigma^*\overline{\sigma^*\overline{L}}= L\,.
\end{equation*}
We note that the holomorphic line bundle
$\sigma^*\overline{\sigma^*\overline{L}}$ is canonically identified
with $L$.

Let $c\, \in\, {\mathbb C}\setminus\{0\}$ be the nonzero
complex number such that
\begin{equation}\label{s.c.c}
(\sigma^*\overline{\alpha})\circ \alpha = c\cdot \text{Id}_L\,,
\end{equation}
where $(\sigma^*\overline{\alpha})\circ \alpha$ is the automorphism
in \eqref{s.c.2}.

We will show that $c$ is actually a real number.

Let
\begin{equation*}
M:= L\otimes \sigma^*\overline{L}
\end{equation*}
be the holomorphic line bundle over $X$. We
note that
\begin{equation*}
\sigma^*\overline{M}=\sigma^*\overline{L}\otimes
\sigma^*\overline{\sigma^*\overline{L}}=
\sigma^*\overline{L} \otimes L = M\,.
\end{equation*}
Let
\begin{equation}\label{s.c.3}
\tau\hbox{:}\ M  \rightarrow
\sigma^*\overline{M}
\end{equation}
be the above isomorphism. It is easy to see that
$(\sigma^*\overline{\tau})\circ\tau$ is the identity
automorphism of the line bundle $M$.

Next we observe that the tensor product of two
homomorphisms $\alpha\otimes \sigma^*\overline{\alpha}$
is an automorphism of the holomorphic line bundle $M$. Indeed,
$\alpha$ sends the line bundle $L$
to $\sigma^*\overline{L}$ and $\sigma^*\overline{\alpha}$ sends
$\sigma^*\overline{L}$ to
$\sigma^*\overline{\sigma^*\overline{L}} = L$. Therefore,
$\alpha\otimes \sigma^*\overline{\alpha}$ sends the line bundle
$L\otimes \sigma^*\overline{L} =: M$ to itself.

Consequently,
\begin{equation}\label{s.c.4}
\alpha\otimes \sigma^*\overline{\alpha} =d\cdot \text{Id}_M\, ,
\end{equation}
where $d$ is a nonzero complex number.\pagebreak

Let
\begin{equation}\label{s.c.5}
\delta := \tau\circ (\alpha\otimes \sigma^*\overline{\alpha})\hbox{:}\ M  \rightarrow \sigma^*\overline{M}
\end{equation}
be the isomorphism, where $\tau$ is the isomorphism in
\eqref{s.c.3} and $\alpha\otimes \sigma^*\overline{\alpha}$
is the automorphism of $M$ in \eqref{s.c.4}.

We consider the automorphism
\begin{equation*}
(\sigma^*\overline{\delta})\circ\delta\hbox{:}\
M \rightarrow M\, ,
\end{equation*}
where $\delta$ is defined in \eqref{s.c.5}. From \eqref{s.c.4}
it follows immediately that
\begin{equation*}
(\sigma^*\overline{\delta})\circ\delta = \vert d\vert^2\cdot
\text{Id}_M\,.
\end{equation*}
On the other hand, from \eqref{s.c.c} it follows that
\begin{equation*}
(\sigma^*\overline{\delta})\circ\delta = c^2\cdot \text{Id}_M\,.
\end{equation*}

Since $c^2\cdot \text{Id}_M=(\sigma^*\overline{\delta})\circ\delta
=  \vert d\vert^2\cdot \text{Id}_M$,
we have $c^2 = \vert d\vert^2$. Consequently,
$c\, \in\, {\mathbb R}$, where $c$ is the scalar in \eqref{s.c.c}.

As $c\, \in\, {\mathbb R}\setminus\{0\}$, the isomorphism
\begin{equation*}
\alpha_0 := \frac{\alpha}{\sqrt{\vert c\vert}}\hbox{:}\ L \rightarrow \sigma^*\overline{L}
\end{equation*}
has the property that
\begin{equation*}
(\sigma^*\overline{\alpha_0})\circ \alpha_0 =
\pm \text{Id}_L\,,
\end{equation*}
where $\alpha$ is the isomorphism in \eqref{s.c.1}.

Therefore, we have proved the following lemma.

\begin{lemma}\label{lem1}
Let $L$ be a holomorphic line bundle over $X$ such that the
holomorphic line bundle $\sigma^*\overline{L}$ is isomorphic
to $L$. Then there is an isomorphism
\begin{equation*}
\alpha\hbox{:}\ L \rightarrow \sigma^*\overline{L}
\end{equation*}
such that $\sigma^*\overline{\alpha}\circ \alpha$ is either
${\rm Id}_L$ or $-{\rm Id}_L$.
\end{lemma}

In the proof of Theorem 4.2 of \cite{AB} it was erroneously
asserted that for any $L$ as in Lemma~\ref{lem1}, there is
always an isomorphism
\begin{equation*}
\alpha\hbox{:}\ L \rightarrow \sigma^*\overline{L}\,.
\end{equation*}
such that $\sigma^*\overline{\alpha}\circ \alpha$ is 
${\rm Id}_L$. In Lemma~\ref{1lemma1} we will see that
this is not the case.

\begin{remark}
{\rm The pair $(X\, ,\sigma)$ corresponds to a geometrically
connected smooth projective curve defined over $\mathbb R$
without any real points. A holomorphic line bundle $L$ over
$X$ of degree $d$ with $\sigma^*\overline{L}$
holomorphically isomorphic to
$L$ corresponds to a real point of the Jacobian $J^d(X)$
of line bundles of degree $d$. If there
is an isomorphism $\alpha\hbox{:}\ L \rightarrow \sigma^*
\overline{L}$ such that $\sigma^*\overline{\alpha}\circ \alpha
= \text{Id}_L$, then $L$ corresponds to an algebraic line
bundle defined over the real algebraic curve.}
\end{remark}

Let ${\mathcal G}$ denote the group of all holomorphic
line bundles $L$ over $X$ such that $\sigma^*\overline{L}$
is holomorphically isomorphic to $L$. The group operation is
given by tensor product. Let
\begin{equation*}
{\mathcal G}_0\, \subset\, {\mathcal G}
\end{equation*}
be the subgroup consisting of all $L\, \in \, {\mathcal G}$
admitting an isomorphism
\begin{equation*}
\alpha\hbox{:}\ L \rightarrow \sigma^*
\overline{L}
\end{equation*}
such that $(\sigma^*\overline{\alpha})\circ \alpha\,=\, \text{Id}_L$.

Using Lemma~\ref{lem1} we have a character
\begin{equation}\label{1chi}
\lambda\hbox{:}\ {\mathcal G} \rightarrow {\mathbb Z}/2{\mathbb Z}
\end{equation}
defined by the following condition:
any $L\, \in \, {\mathcal G}$ admits an isomorphism
\begin{equation*}
\alpha\hbox{:}\ L \rightarrow \sigma^*
\overline{L}
\end{equation*}
such that $(\sigma^*\overline{\alpha})\circ \alpha\,=\, \lambda(L)\cdot
\text{Id}_L$ after identifying ${\mathbb Z}/2{\mathbb Z}$
with $\pm 1$. Note that for any nonzero complex number $z$
and any isomorphism $\alpha$ as above,
we have $(\sigma^*\overline{z\alpha})\circ (z\alpha) =
\vert z\vert^2\cdot (\sigma^*\overline{\alpha})\circ \alpha$.
Using this it follows immediately that the above map $\lambda$ is
well-defined.

The kernel of the homomorphism $\lambda$ in \eqref{1chi} coincides
with the subgroup ${\mathcal G}_0$.

\begin{lemma}\label{1lemma1}
The homomorphism $\lambda$ defined in \eqref{1chi} is surjective.
\end{lemma}

\begin{proof}
Since the anti-holomorphic involution $\sigma$ of $X$ does not
have any fixed points, there exists a meromorphic function $f$
on $X$ such that
\begin{equation}\label{1idf}
f\cdot \overline{f\circ\sigma} = -1
\end{equation}
(see \cite{Wi} for the construction of $f$). We note that
$\overline{f\circ\sigma}$ is also a holomorphic function over $X$.

From \eqref{1idf} it follows that $f$ is a nonconstant function.
Let $D_0$ (respectively, $D_1$) be the effective
divisor on $X$ defined
by the zeros (respectively, poles) of $f$. From \eqref{1idf}
it follows that
\begin{equation}\label{1f}
\sigma(D_0)  = D_1\,.
\end{equation}

Let $L = {\mathcal O}_X(D_0)$ be the holomorphic line bundle
over $X$ defined by the divisor $D_0$. Similarly, let
$L' = {\mathcal O}_X(D_1)$ be the holomorphic line bundle
defined by $D_1$. From \eqref{1f} it follows that
$\sigma^* \overline{L} = L'$.

The meromorphic function $f$ defines an isomorphism of
${\mathcal O}_X(D_0)$ with ${\mathcal O}_X(D_1)$. Let
$\alpha$ be the isomorphism of $L$ with
$\sigma^* \overline{L} =
L'$ given by $f$. From \eqref{1idf} it follows that
$(9\sigma^*\overline{\alpha})\circ \alpha\,=\, -\text{Id}_L$.
Therefore, the character $\lambda$ in \eqref{1chi} is nontrivial.
This completes the proof of the\break lemma.\hfill $\Box$
\end{proof}

The following lemma gives the parity of the degree
of any line bundle $L\,\in\,{\mathcal G}$ over $X$ with
$\lambda(L)  = -1$.

\begin{lemma}\label{1lem2}
Let $L\, \in\, {\mathcal G}$ be a holomorphic line bundle
over $X$ such that $L\, \notin\, {\mathcal G}_0$, i.e.{\rm ,}
$\lambda(L)\,=\, -1$. Then
${\rm degree}(L) \, \equiv\, {\rm genus}(X)+ 1$ {\rm mod} $2$.
\end{lemma}

\begin{proof}
Replacing $L$ by the tensor product of
$L$ with the holomorphic line bundle ${\mathcal O}_X(D)$,
where $D$ is an effective divisor on $X$
such that $\sigma(D) = D$
and $\text{degree}(D)$ is sufficiently large,
we may assume that $H^1(X,\, L) = 0$. Note that
since $\sigma$ is fixed-point free, the degree of $D$ is even.

Fix an isomorphism $\alpha\hbox{:}\ L \rightarrow \sigma^*
\overline{L}$
such that $\sigma^*\overline{\alpha}\circ \alpha\,=\, -\text{Id}_L$.
This isomorphism $\alpha$ induces a conjugate linear automorphism 
$\widehat{\alpha}$ of the complex vector space $H^0(X,\, L)$ such that 
$\widehat{\alpha}^2 = -\text{Id}$. From this it follows
immediately that the complex dimension of $H^0(X,\, L)$ is even.

The Riemann--Roch theorem says
\begin{equation*}
\dim H^0(X,\, L) -\dim H^1(X,\, L) = \text{degree}(L)
-\text{genus}(X) +1\,.
\end{equation*}
Since $\dim H^0(X,\, L)$ is even and $H^1(X,\, L) \,=\, 0$,
the lemma follows from the Riemann--Roch theorem.\hfill $\Box$
\end{proof}

\begin{remark}
{\rm Take any line bundle $L\, \in\, {\mathcal G}\backslash{\mathcal G}_0$. From
Lemma~\ref{1lem2} it follows that
${\rm degree}(L) = 2d_0$ if the genus of $X$ is odd,
and ${\rm degree}(L) = 2d_0 + 1$ if the genus is even.
Take a divisor $D_0 =\sum_{i=1}^{d_0}x_i X$ of degree $d_0$ such that $\sigma (D_0) \,=\, D_0$ and
$\{x_1, \dots, x_{d_0}\}$ are distinct points.
Consider the line bundle $L':= L\otimes {\mathcal O}_X(-D_0)$
over $X$. Note that $L'\, \in\, {\mathcal G}\setminus
{\mathcal G}_0$, and ${\rm degree}\ (L')$ is zero or one depending
on the parity of ${\rm genus}\ (X)$.}
\end{remark}

Consider the intersection $J^0(X)\cap {\mathcal G}_0$.
The homomorphism $\rho$ in p.~147, eq.~(3) of \cite{AB} maps
to it. Combining Lemmas~\ref{1lemma1} and \ref{1lem2}
we have the following corrected version of Theorem 4.2 of \cite{AB}.

\begin{thm}[\!]\label{thm1}
The image of the homomorphism $\rho$ in p.~$147${\rm ,} eq.~$(3)$ of {\rm \cite{AB}}
coincides with $J^0(X)\cap {\mathcal G}_0$. If the genus of $X$
is even{\rm ,} then the image of the homomorphism $\rho$ in
p.~$147${\rm ,} eq.~$(3)$ of {\rm \cite{AB}} coincides with $J^0(X)\cap {\mathcal G}$.
If the genus of $X$ is odd{\rm ,} then the image of $\rho$ is a subgroup
of $J^0(X)\cap {\mathcal G}$ of index two.
\end{thm}

\section{The Jacobian as dual space of holomorphic forms}

Following the setting and notation of \cite{AB}, let $J_1(X)$ denote the
Jacobian of $X$ obtained as the quotient of $H^0(X,\Omega)^*$ (the dual
of the space of holomorphic forms on $X$) by the action of
$H_1(X,\mathbb Z)$.
In Proposition~3.1 of \cite{AB} we showed that if $\{\gamma_1,\dots,\gamma_g,
\delta_1,\dots,\delta_g\}$ is a canonical (symplectic) basis of
$H_1(X,\mathbb Z)$ satisfying
\begin{equation}\label{invhomology}
\sigma_\#(\gamma_j) = \gamma_j,
\end{equation}
then the associated basis of holomorphic form, $\{\omega_1,
\dots, \omega_g\}$ ($\int_{\gamma_j} \omega_k = \delta_{jk}$)
is invariant, that is, $\overline{\sigma^*(\omega_j)} =
\omega_j$. Here $\sigma_\#$ and $\sigma^*$ denote the maps induced by
$\sigma$ in homology and 1-forms respectively.

Let $\pi\hbox{:}\ \mathbb C^g \to J_1(X) \cong \mathbb C^g / \mathbb Z^g$ be the
natural projection, and let $\sigma_1$ denote the involution induced by
$\sigma$ on $J_1(X)$. In \cite{AB} we showed that the lift of $\sigma_1$
to $\mathbb C^g$, with the above basis, is given by conjugation. The
fixed points of $\sigma_1$ are therefore given by the solutions of the
equation
\begin{equation}\label{fixedpoints}
\overline{z} = z + n + Pm,
\end{equation}
where $P$ is the period matrix and $n$ and $m$ are points in
$\mathbb Z^g$. We obtain the solutions to this equation in two
different ways, depending on whether the genus of $X$ is even or odd.

Let us first fix some notation. Let $I_n$ denote the $n \times n$
identity matrix, and $K_n$ the $n \times n$ matrix with entries equal
to $1$ in the anti-diagonal and $0$ in all other entries.

Assume first that $X$ has even genus. The fundamental group of the Klein
surface $Y$ has one relation given by $c^2(a_1b_1a_1^{-1}b_1^{-1})
\cdots (a_{g-1}b_{g-1}a_{g-1}^{-1}b_{g-1}^{-1}) = e$, from which we get,
by a simple topological argument (see \cite{Bl}), a symplectic basis of
$X$, say $\mathcal B$. It is easy to see that the action of $\sigma_\#$
on $H_1(X,\mathbb Z)$ with respect to $\mathcal B$ is given by the
matrix $K_{2g}$. We make a change of basis using the $2g \times 2g$ matrix
\begin{equation*}
C = \begin{pmatrix} -I_g & (I+K)_g\\[.3pc]
-K_g & K_g\end{pmatrix};
\end{equation*}
let $\mathcal B_1 = \mathcal B C = \{\gamma_1, \dots, \gamma_g,
\delta_1, \dots, \delta_g \}$ be the new basis.
Since $C$ satisfies $C^tJC = J$, where $J = \left(\begin{smallmatrix}0 &
-I_g\\ I_g & 0\end{smallmatrix}\right)$ is the standard intersection
matrix, we have that $\mathcal B_1$ is symplectic.
The action of $\sigma_\#$ with respect to $\mathcal B_1$ is given by
\begin{equation*}
C^{-1}\sigma_\#C = \begin{pmatrix}I_g & (-2I-K)_g \\[.3pc]
0 & -I_g\end{pmatrix},
\end{equation*}
so $\mathcal B_1$ satisfies condition \eqref{invhomology}.

Let $P$ denote the period matrix with respect to this new basis of
$H_1(X,\mathbb Z)$ and the associated basis of holomorphic forms;
denote by $p_{jk}$ its entries. Let $A$ be the $g \times g$ matrix given
by $A = -2I_g - K_g$ and denote by $a_{jk}$ the entries of $A$. Then
we have
\begin{equation*}
\overline{p_{kj}} = \overline{\int_{\delta_j} \omega_k} =
\int_{\delta_j} \overline{\omega_k} = \int_{\delta_j} \sigma^*(\omega_k)
= \int_{\sigma_\#(\delta_j)} \omega_k = a_{jk} - p_{jk}.
\end{equation*}
So the real part of $P$ is equal to $\hbox{Re}\,P = -I - \frac{1}{2}K$.

Since the imaginary part of $P$ is invertible (see for example
Proposition III.2.8 of \cite{FK}), we have that any point $z$ of
$\mathbb C^g$ can be written as $z = x + P y$, where $x$ and $y$ are
points in $\mathbb R^g$. Considering the real and imaginary parts
of \eqref{fixedpoints} we obtain the following two equations:
\begin{equation*}
\begin{cases} 0 & = n + (\hbox{Re}\,P)m,\\[.2pc]
-(\hbox{Im}\,P) y  & = (\hbox{Im}\,P)y + (\hbox{Im}\,P)m.
\end{cases}
\end{equation*}
From the first equation, using the expression of $\hbox{Re}\,P$ obtained above
and the fact that all entries of $n$ are integers we get that the
entries of $m$ are even integers. Since $\hbox{Im}\,P$ is invertible, the
second equation gives $y = -\frac{1}{2}m$, which implies that $y$ has
integer entries. The set of fixed points of $\sigma_1$ is therefore
given by the projection of
\begin{equation*}
\mathcal T = \{x + Py;~x \in \mathbb R^g,~y\in \mathbb Z^g\} \subset
\mathbb C^g
\end{equation*}
to $J_1(X)$.

The odd genus case is handled in a similar way (although the
computations are a little more complicated): using topological
arguments we get a symplectic basis $\mathcal B$ of $H_1(X,\mathbb Z)$;
the action of $\sigma_\#$ with respect to this basis is given by
\begin{equation*}
\begin{pmatrix}1 & 0 & 0 & 0\\
0 & 0 & 0 & K_{g-1}\\
0 & 0 & -1 & 0\\
0 & K_{g-1} & 0 & 0\end{pmatrix}.
\end{equation*}
We change basis in $H_1(X,\mathbb Z)$ to $\mathcal B_1 = \mathcal B C$,
where
\begin{equation*}
C = \begin{pmatrix}1 & 0 & 0 & 0\\
0 & -I_{g-1} & 0 & (I+K)_{g-1}\\
0 & 0 & 1 & 0\\
0 & -K_{g-1} & 0 & K_{g-1}\end{pmatrix}.
\end{equation*}
One can easily show that $\mathbb B_1$ is symplectic and satisfies 
condition \eqref{invhomology}. The period matrix satisfies the identity
$P = A - \overline{P}$, where $A$ is the matrix
$A = \left( \begin{smallmatrix}0 & 0\\ 0 & (-2I-K)_{g-1}\end{smallmatrix}
\right)$. Splitting \eqref{fixedpoints} in its real and imaginary parts
as above we obtain the following pairs of equations (we used again the
fact that $\hbox{Im}\,P$ is invertible):
\begin{equation*}
\begin{cases} n & = -(\hbox{Re}\,P)m,\\
y & = -\frac{1}{2} m.\end{cases}
\end{equation*}
From the first equation we get that if $m^t = (m_1, \dots, m_g)
\in \mathbb Z^g$ then $m_2, \dots, m_g$ are even integers. So the set of
fixed points of $\sigma_1$ has two components, given by the projections
of the following two sets to $J_1(X)$:
\begin{align*}
\mathcal T_1 & = \{x + Py;~x \in \mathbb R^g,~y\in \mathbb Z^g\} \subset
\mathbb C^g,\\[.5pc]
\mathcal T_2 & = T_1 + P\,\left(\frac{1}{2},0, \dots, 0\right)^t =
\mathcal T_1 +
\frac{1}{2}p_1 \subset \mathbb C^g.
\end{align*}
Here $p_1$ is the first period. We refer the reader to \cite{NT} where
similar results are stated, although with different computations. 
The Jacobian of $Y$ can then be identified with the subgroup 
$\pi(\mathcal T_1)$ of $J_1(X)$.

\section{Divisors}

Let $J_0(X)$ be the Jacobian of $X$ given as the divisor class group,
that is degree zero divisors quotiented by the principal divisors. We
have a natural involution in $J_0(X)$ defined by $\sigma_0([D]) =
[\sigma^*(D)]$, where square brackets denote equivalence classes and
$\sigma^*$ the natural extension of $\sigma$ to divisors.
It was shown in \cite{AB} that the involution $\sigma_0$ is equivalent
to $\sigma_1$ by~the Abel--Jacobi map. A divisor class $[D]$ is fixed by
$\sigma_0$ if $D$ is linearly equivalent to $\sigma^*(D)$; that is,
there exists a holomorphic function $f\hbox{:}\ X \to \widehat{\mathbb C}$
such that $D - \sigma^*(D) = \hbox{div}(f)$. This implies that
$\hbox{div}(\overline{f\circ \sigma}) = -\hbox{div}(f)$, so there exists a constant
$c$ such that $f\cdot\overline{f\circ\sigma} = c$ (here $\cdot$ denotes
multiplication of complex numbers). But then $\overline{f\circ\sigma}
\cdot f = \overline{c} = c$, so $c$ is a real number. Multiplying $f$ by
a number we can assume that $c = \pm 1$.

If $c = 1$, let $h$ be the function $h = f+1$. Then
$\overline{h\circ\sigma} = h/f$ so that $f = h / (\overline{h\circ\sigma})$,
which implies $\hbox{div}(f) = \hbox{div}(h) - \hbox{div}(\overline{h\circ\sigma}) =
\hbox{div}(h) - \sigma^*(\hbox{div}(h))$. Let $E$ be the divisor $E = D - \hbox{div}(h)$.
Then we have that $E$ is linearly equivalent to $D$ and
\begin{align*}
\sigma^*(E) = \sigma^*(D) - \sigma^*(\hbox{div}(h)) = (D - \hbox{div}(f)) - (\hbox{div}(h) -
\hbox{div}(f)) = E.
\end{align*}
Therefore $D$ is linearly equivalent to a divisor ($E$) that comes from
the surface $Y$.

If $f\cdot \overline{f\circ\sigma} = -1$ it can be easily seen that $D$
is not linearly equivalent to any divisor that comes from $Y$. So we get that the
set of fixed points on $\sigma_0$ consists of the following disjoint sets:\pagebreak
\begin{align*}
\mathcal T_1 & = \{[D];~\deg(D) = 0,~\sigma^*(D) = D\}~\mathrm{and}~\\[.4pc]
\mathcal T_2 & = \{[D];~\deg(D) = 0,~D-\sigma^*(D) = \hbox{div}(f),~f\cdot
\overline{f\circ\sigma} = 1\}.
\end{align*}
Clearly $\mathcal T_1$ is not empty.

If $X$ has even degree, then by Lemma~\ref{1lemma1} we have that
$\mathcal T_2$ is empty. In the case of odd degree, let $f$ be as in
\eqref{1idf} and let $D_0$
and $D_1$ denote the divisors defined by the zeroes and poles of $f$,
respectively. Equation \eqref{1idf} gives $\sigma^*(D_1) = D_0$.
Again by Lemma~\ref{1lemma1} we have that $D_0$ has even degree.
Let $E$ be a divisor with $\hbox{degree}(E) = \hbox{degree}(D_0)$ and $E = \sigma^*(E)$
(observe that this last condition, since $\sigma$ does not have fixed
points, forces the degree of $E$ to be even). Define $\mathcal X =
E - D_0$. Then $\mathcal X$ is a degree zero divisor satisfying
$\mathcal X - \sigma^*(\mathcal X) = \hbox{div}(f)$, so $[\mathcal X] \in
\mathcal T_2$. This shows that $\mathcal T_2$ is not empty in the case
of $X$ having odd degree. It is easy to see that $[D] \in \mathcal T_1$
if and only if $[D + \mathcal X] \in \mathcal T_2$, so $\mathcal T_2$ is
the translation of $\mathcal T_1$ by $\mathcal X$.

We can identify $\mathcal T_1$ with the Jacobian of $Y$ (and
with the component $\mathcal T_1$ obtained in the previous version, via
the Abel--Jacobi mapping).

\end{document}